\documentclass[11pt]{article}

\usepackage{amsthm,amsmath,natbib, amssymb}
\RequirePackage[colorlinks,citecolor=blue,urlcolor=blue]{hyperref}

\RequirePackage{hypernat}

\font\tenmath=msbm10
\font\sevenmath=msbm7
\font\fivemath=msbm5
\newfam\mathfam \textfont\mathfam=\tenmath
\scriptfont\mathfam=\sevenmath \scriptscriptfont\mathfam=\fivemath

\usepackage{amsmath, moreverb, hyperref}
 \hypersetup{colorlinks=true,citecolor=blue}

\usepackage{natbib, color}
\usepackage{amsfonts, fancybox}
\usepackage{amssymb}
\usepackage[american]{babel}
\usepackage{graphicx, pstricks, pst-node, color}
\usepackage{psfrag}\usepackage{subfigure}
\usepackage{flafter}
\usepackage[section]{placeins}

\newcommand{\cB}{\mathcal{B}}
\newcommand{\cC}{\mathcal{C}}

\newcommand{\cF}{\mathcal{F}}
\newcommand{\cG}{\mathcal{G}}

\newcommand{\cN}{\mathcal{N}}

\newcommand{\Pro}{\mathbf{P}\kern-0.12em}

\newcommand{\R}{{\rm I}\kern-0.18em{\rm R}}
\newcommand{\h}{{\rm I}\kern-0.18em{\rm H}}
\newcommand{\K}{{\rm I}\kern-0.18em{\rm K}}
\newcommand{\p}{{\rm I}\kern-0.18em{\rm P}}
\newcommand{\E}{{\rm I}\kern-0.18em{\rm E}}
\newcommand{\Z}{{\rm Z}\kern-0.18em{\rm Z}}

\newcommand{\KL}{\mathsf{\mathop{KL}}}

\newcommand{\tr}{\mathbf{Tr}}

\newcommand{\pn}{\p_{\kern-0.25em n}}
\newcommand{\pnm}{\p_{\kern-0.25em n,m}}
\newcommand{\psubm}{\p_{\kern-0.25em m}}
\newcommand{\psubp}{\p_{\kern-0.25em p}}
\newcommand{\cfi}{\cF_{\kern-0.25em \infty}}

\newcommand{\card}{\mathop{\mathrm{card}}}
\newcommand{\argmin}{\mathop{\mathrm{argmin}}}

\newcommand{\sign}{\mathop{\mathrm{sign}}}
\newcommand{\ud}{\mathrm{d}}

\newcommand{\epr}{\hfill\hbox{\hskip 4pt\vrule width 5pt
                  height 6pt depth 1.5pt}\vspace{0.5cm}\par}
\newcommand{\eps}{\varepsilon}

\newtheorem{TH1}{Theorem}

\newtheorem{cor}{Corollary}
\newtheorem{lem}{Lemma}

\newlength{\minipagewidth}
\graphicspath{{./eps/}}
\begin{document}
\title{Estimation of Covariance Matrices under Sparsity Constraints}
\author{\normalsize
Philippe Rigollet \thanks{Department of Operations Research and Financial Engineering, Princeton University, Princeton, NJ 08544, USA. Supported in part by NSF grants DMS-0906424 and CAREER-DMS-1053987. }
\and \normalsize Alexandre B. Tsybakov \thanks{
Laboratoire de Statistique, CREST-ENSAE,
 3, av. Pierre Larousse, F-92240 Malakoff Cedex, France. Supported in part by ANR under the grant ``Parcimonie".}}
\date{\normalsize \today}
\maketitle


\setcounter{equation}{0}
\label{SEC:intro}
{\bf Introduction.} Estimation of covariance matrices in various norms is a critical issue that finds applications in a wide range of statistical problems, and especially in principal component analysis. It is well known that, without further assumptions, the empirical covariance matrix $\Sigma^*$ is the best possible estimator in many ways, and in particular in a minimax sense. However, it is also well known that $\Sigma^*$ is not an accurate estimator when the dimension $p$ of the observations is high.  The minimax analysis carried out by Tony Cai and Harry Zhou ([CZ] in what follows)  guarantees that for several classes of matrices with reasonable structure (sparse or banded matrices), the fully data-driven thresholding estimator achieves the best possible rates when $p$ is much larger than the sample size $n$. This is done, in particular, by proving minimax lower bounds that ensure that no estimator can perform better than the hard thresholding estimator, uniformly over the sparsity classes $\cG_q$ for each $0\le q<1$. This result has a flavor of universality in the sense that one and the same estimator is minimax optimal for several classes of matrices.  

Our comments focus on the sparsity classes of matrices.
\begin{itemize}
\item[($a$)] {\it Optimal rates.} Optimal rates are obtained in [CZ] under the assumption that the dimension is very high: $p\ge n^{\nu}$, $\nu>1$. Thus, the
case of dimensions smaller than $n$, or even $p\approx n$, is excluded. This seems to be due to the technique of proving the lower bound (Theorem 2 in [CZ]). Indeed, by a different technique, we show that the lower bound holds without this assumption, cf. Theorem 1 below. Furthermore, in general, our lower rate $\psi^{(1)}$ is different from that obtained in [CZ] and has ingredients similar to the optimal rate for the Gaussian sequence model. We conjecture that it is optimal for all admissible configurations of $n,p$, and sparsity parameters.   
\item[($b$)] {\it Frobenius norm and global sparsity.} We argue that the Frobenius norm is naturally adapted to the structure of the problem, at least for Gaussian observations, and we derive optimal rates under the Frobenius risk and global sparsity assumption. 
\item[($c$)] {\it Approximate sparsity.} Again under the Frobenius risk, one can obtain not only the minimax results but also oracle inequalities. We demonstrate it for the soft-thresholding estimator. This allows us to deal with a more general setup where the covariance matrix is not necessarily sparse but can be well approximated by a sparse matrix.
\end{itemize}


Below we denote by $\|A\|$  the Frobenius norm of a matrix $A$:
$$
\|A\|^2=\tr(AA^\top)=\sum_{i,j}a_{ij}^2\,
$$
where $\tr(B)$ stands for the trace of square matrix $B$.
Moreover, for $q>0$, we denote by $|v|_q$ the $\ell_q$-norm of a vector $v$ and by $|A|_q$ the $\ell_q$ norm of the off-diagonal entries of $A$. We set $|A|_0=\sum_{i\ne j} I(a_{ij}\ne 0)$ (the number of non-zero off-diagonal entries of $A$). The operator $\ell_q \to \ell_q$ norm of $A$ is denoted by $\|A\|_q$.


\medskip

{\bf Frobenius norm and Sparsity.} The cone of positive semi-definite (PSD) matrices can be equipped with a variety of norms, even more so than a vector space. [CZ] choose the $\|\cdot\|_1$ norm and consider classes of matrices that are essentially adapted to this metric. For example, the class $\cG_q$ defined in (1) controls the largest $\ell_q$ norm of the columns of the covariance matrix $\Sigma$ with $0\le q<1$ while the $\|\cdot\|_1$ norm measures the largest $\ell_1$ norm of the columns of $\hat \Sigma -\Sigma$. Theorem~\ref{TH:low} below indicates that for $q =1$ consistent estimators do not exist. 

One may wonder whether faster rates can be obtained if, for example, $\Sigma$ has one row/column with large $\ell_q$ norm and all other rows/columns have small $\ell_q$ norm. It is quite clear that the $\|\cdot\|_1$ norm fails to capture such a behavior and we need to resort to other norms. As we see below, this is achievable when the Frobenius norm is used.

The Frobenius norm is a rather weak norm on the PSD cone. Indeed, it is very much a vector norm unlike the $\|\cdot\|_1$ norm used by [CZ] or the spectral norm that are operator norms. However, the choice of a norm is rather subjective but some general guidelines exist in a given statistical setup. It can be motivated by the idea of minimizing the Kullback-Leibler divergence between the true distribution and its estimator~\citep[see, e.g.,][]{Rig12}. This principle naturally gives rise to the use of the Frobenius norm in Gaussian covariance matrix estimation, as indicated by the following lemma.
\begin{lem}
\label{LEM:KL}
Let  $I_p$ be the $p \times p$ identity matrix and $\Delta$ be a symmetric $p \times p$ matrix such that $I_p +  \Delta$ is PSD. Denote by $P_\Sigma$ the distribution of $\cN_p(0, \Sigma)$ (a zero-mean normal random variable in $\R^p$ with covariance matrix $\Sigma>0$). Then, for any $0<\eps<1$, the Kullback-Leibler divergence between $P_{I_p + \eps\Delta}$ and $P_{I_p}$ satisfies
$$
\KL(P_{I_p +  \eps \Delta}, P_{I_p}) \le \frac{g(-\eps)}{2}\|\Delta\|^2\,,
$$
where 
$$
g(\eps)=\frac{\eps - \log(1+\eps)}{\eps^2}\,.
$$
Moreover if $\|\Delta\|_2 \le 1$, we have
\begin{equation}\label{kull}
\KL(P_{I_p + \eps \Delta}, P_{I_p}) \ge \frac{(1-\log 2)\eps^2}{2} \|\Delta\|^2\,.
\end{equation}
\end{lem}
{\sc Proof.} 
Take $\Sigma= I_p + \eps \Delta$ and observe that 
$$
\KL(P_\Sigma, P_{I_p})=\E \log\left(\frac{\ud P_\Sigma}{\ud P_{I_p}}(X)\right)=\frac{1}{2}\E \log \left( \frac{1}{\det(\Sigma)}\right) + \frac{1}{2}\E[X^\top X-X^\top \Sigma^{-1} X]\,,
$$
where $X \sim \cN_p(0, \Sigma)$. Let $\lambda_1, \ldots, \lambda_p$ denote the eigenvalues of $\Delta$ and recall that $\det(\Sigma)=\prod_j(1+\eps \lambda_j)$. Moreover,
$$
\E[X^\top X-X^\top \Sigma^{-1} X]=\tr(\E[ X X^\top] - \Sigma^{-1} \E[ X X^\top])=\tr(\Sigma-I_p)= \sum_j \eps \lambda_j\,.
$$
Therefore, 
$$
\KL(P_\Sigma, P_{I_p})=\frac{1}{2}\sum_{j=1}^p[\eps \lambda_j-\log ( 1+\eps \lambda_j) ]\le \frac{1}{2}\sum_{j=1}^p g(\eps\lambda_j)\lambda_j^2\,.
$$
Note now that since $I_p+\Delta$ is PSD, then $\lambda_j \ge -1$ for all $j =1, \ldots, p$. Therefore, since $g$ is monotone decreasing on $(-1, \infty)$, it yields $g(\eps \lambda_j) \le g(-\eps)$. The second statement of the lemma follows by observing that if $\|\Delta\|_2 \le 1$, then $\eps \lambda_j \le \eps \le 1$ for all $j=1, \ldots, p$\,.
\epr

{\bf Minimax lower bounds over classes of sparse matrices.} We denote by $\sigma_{ij}$ the elements of $ \Sigma$ and by $\sigma_{(j)}$ the $j$th column of $ \Sigma$ with its $j$th component replaced by 0. For any $q> 0, R>0$, we define the following classes of matrices: 
$$
\cG_q^{(0)}(R)=\left\{ \Sigma\in  \cC_{>0}\,:\, |\Sigma|_q^q \le R\,,\   \sigma_{ii} = 1, \forall i\right\}\,,$$
$$
\cG_q^{(1)}(R)=\left\{ \Sigma\in  \cC_{>0}\,:\, \max_{1\le j\le p} |\sigma_{(j)}|_q^q \le R\,, \  \sigma_{ii} = 1, \forall i\right\}\,,
$$
where $\cC_{>0}$ is the set of all positive definite symmetric $p\times p$ matrices. For $q=0$, we define the classes $\cG_0^{(0)}(R)$ and $\cG_0^{(1)}(R)$ analogously, with the respective constraints $|\Sigma|_0 \le R$ and $\max_{1\le j\le p} |\sigma_{(j)}|_0 \le R$. Here $R$ is an integer for the class $\cG_0^{(1)}(R)$, and an even integer for $\cG_0^{(0)}(R)$ in view of the symmetry. We assume that $R=2k\le {p(p-1)}$ for $\cG_0^{(0)}(R)$ and $R=k\le {p-1}$ for $\cG_0^{(1)}(R)$ where $k$ is an integer. Set
$$
\psi^{(0)}=R^{1/2}\left(\frac{1}{n}\log\left(1+\frac{c_0p^2}{Rn^{q/2}}\right)\right)^{\frac1{2}-\frac{q}{4}}, \quad  \psi^{(1)}
=R\left(\frac{1}{n}\log\left(1+\frac{c_0p}{Rn^{q/2}}\right)\right)^{\frac{1-q}{2}}\,,
$$
for some positive constant $c_0$ that does not depend on the parameters $p, n, R$.
The following minimax lower bounds hold.
\begin{TH1}
\label{TH:low}
Fix $R>0$, $0 \le q  \le 2$, $C_0>0$ and integers $n\ge 1$, $p\ge 2$. Consider the conditions
\begin{eqnarray}\label{EQ:low:ass}
&&R((\log p)/n)^{1-q/2}\le C_0, \quad \quad R((\log p)/n)^{(1-q)/2}\le C_0, \quad \quad  R^{-1}((\log p)/n)^{q/2}\le C_0.
\end{eqnarray}  
Let $X_1, \ldots, X_n$ be i.i.d.~$\cN_p(0, \Sigma)$ random vectors, and let $w:[0,\infty)\to [0,\infty)$ be a monotone non-decreasing function such that $w(0)=0$ and $w\not\equiv 0$. Then there exist constants $c_0>0,c_1>0, c>0$ depending only on $C_0$ such that, under the first and third conditions in~\eqref{EQ:low:ass}, 
\begin{eqnarray}
\label{EQ:low1}
&&\inf_{\hat \Sigma} \sup_{\Sigma \in \cG_q^{(0)}(R)}E_{\Sigma} w\big( \big\| \hat \Sigma - \Sigma \big\|/c_1\psi^{(0)})\ge c,
\end{eqnarray}
and under the second and third conditions in~\eqref{EQ:low:ass}, 
\begin{eqnarray}
\label{EQ:low2}
&&\inf_{\hat \Sigma} \sup_{\Sigma \in \cG_q^{(1)}(R)}E_{\Sigma} w\big( \big\| \hat \Sigma - \Sigma \big\|_1/c_1 \psi^{(1)})\ge c, \quad \forall \ 0\le q\le 1,
\end{eqnarray}
where $E_{\Sigma}$ denotes the expectation with respect to the joint distribution of $X_1, \ldots, X_n$ and the infimum is taken over all estimators   based on $X_1, \ldots, X_n$. 
\end{TH1}

\noindent{\sc Proof.} 
We first prove \eqref{EQ:low1} with $q=0$ and $R=2k$. Assume first that $k\le p^2/16$. We use Theorem~2.7 in~\citet{Tsy09}. It is enough to check that there exists a finite subset $\cN$ of $\cG_0^{(0)}(2k)$ such that, for some constant $C>0$ and some $ \psi\ge  C\psi^{(0)}$, we have
$$
\begin{array}{rll}
(i) & \|\Sigma - \Sigma'\|\ge \psi\,, & \forall \ \Sigma \neq  \Sigma' \in \cN\cup \{I_p\},\\
(ii) & n \,\KL(P_\Sigma, P_{I_p}) \le  2^{-4}\log(\card \cN)\,, & \forall \ \Sigma \in \cN\,.
\end{array}
$$
We show that these conditions hold for 
$$
\psi=\left(\frac{k}{n}\log\left(1+\frac{ep(p-1)}{2k}\right)\right)^{1/2}\,.
$$
Let $\cB$ be the family of all $p \times p$ symmetric binary matrices, banded such that for all $B \in \cB$,  $b_{ij}=0$ if $|i-j|>\sqrt{k}$, with $0$ on the diagonal and exactly $k$ nonzero over-diagonal entries equal to 1. Let $M$ be the number of elements in the over-diagonal band where the entry $1$ can only appear. For $k \le p^2/4$ we have $M\ge p\sqrt{k}-k \ge p\sqrt{k}/2$. Therefore for, $k\le p\sqrt{k}/4$, Lemma~A.3 in~\citet{RigTsy11} implies that there exists a subset $\cB_0$ of $\cB$ such that  for any $B, B' \in \cB_0, B\neq B'$, we have $\|B-B'\|^2 \ge (k+1)/4$,   and
\begin{equation}
\label{EQ:log_card}
\log(\card \cB_0) \ge C_1k\log\left(1+\frac{ep}{4\sqrt{k}}\right)
\end{equation}
for some absolute constant $C_1>0$. Consider the family of matrices
$
\cN=\{
 \Sigma=I_p+\frac{a}{2}B: \ B \in \cB_0\}\,,
$
where 
$$
a=a_0\left(\frac{1}{n}\log\left(1+\frac{ep}{4\sqrt{k}}\right)\right)^{1/2}\,
$$
for some $a_0>0$. All matrices in $\cN$ have at most $2\sqrt{k}$ nonzero elements equal to $a$ in each row. Therefore, the first inequality in~\eqref{EQ:low:ass} guarantees that for $a_0$ small enough, matrices $I_p+aB$ with $B \in \cB_0$ and, a fortiori, $\Sigma \in \cN$ are diagonally dominant and hence PSD.  Thus, $\cN \subset \cG_0^{(0)}(2k)$ for sufficiently small $a_0>0$.  Also, for any $\Sigma, \Sigma' \in \cN$, $\Sigma \neq \Sigma'$, we have
$$
\|\Sigma - \Sigma'\|^2 \ge C_2a^2k
$$
for some absolute constant $C_2>0$. It is easy to see that this inequality also holds with a different $C_2$ if $\Sigma$ or $\Sigma'$ is equal to $I_p$. The above display implies $(i)$. 
To check $(ii)$, observe first that since $I_p+aB$ is PSD, we can apply Lemma~\ref{LEM:KL} with $\Delta=aB$, $\eps=1/2$, to get
$$
n\KL(P_\Sigma, P_{I_p})\le \frac{na^2g(-1/2)}{2}\|B\|^2\le a^2kn\,,\quad \forall \ \Sigma\in \cN\,.
$$
To prove $(ii)$, it suffices to take $a_0^2< 2^{-4}C_1$, and to use \eqref{EQ:log_card}.  This proves \eqref{EQ:low1} with $q=0$ under the assumption $k\le p^2/16$. 
  The case $q=0$, $k> p^2/16$ corresponds to a rate $\psi^{(0)}$ of  order $\sqrt{p/n}$ and is easily treated via the Varshamov-Gilbert argument (we omit the details). 

Next, observe that~\eqref{EQ:low1} for $0<q\le 2$, follows from the case $q=0$. Indeed, let $k$ be the maximal integer such that $2ka^q \le R$ (we assume $a_0$  small enough to have $k\ge 1$, cf. the third inequality in~\eqref{EQ:low:ass}). Hence,  $|\Sigma|_q^q=2ka^q\le R$ for any $\Sigma \in \cN$. Also, $a\sqrt{k} \le R^{1/2}a^{1-q/2}/\sqrt{2}$ and thus the first inequality in~\eqref{EQ:low:ass} ensures the positive definiteness of all $\Sigma \in \cN$ for small $a_0$.
For this choice of $k$, we have $k+1>Ra^{-q}/2$ and $k \le C_3R n^{q/2}$ with some constant $C_3>0$. It can be easily shown that $(i)$ holds with 
$$
\psi^2\ge CR\left(\frac{1}{n}\log\left(1+\frac{ep^2}{Ra^{-q}}\right)\right)^{1-\frac{q}{2}}\ge C R\left(\frac{1}{n}\log\left(1+\frac{c_0p^2}{Rn^{q/2}}\right)\right)^{1-\frac{q}{2}}.
$$

The proof of~\eqref{EQ:low2} is quite analogous, with the only difference that $\cB$ is now defined as the set of all symmetric binary matrices  with exactly $k$ off-diagonal entries equal to 1 in the first row and in the first column and all other entries 0. Then, for $k\le (p-1)/2$, Lemma~A.3 in~\citet{RigTsy11} implies that there exists a subset $\cB_1$ of $\cB$ such that  for any two distinct $B, B' \in \cB_1$, we have $|b_{(1)}-b'_{(1)}|_1\ge (k+1)/4$   (consequently, $\|B-B'\|_1\ge (k+1)/4$) and
\begin{equation}
\label{EQ:log_card1}
\log(\card \cB_1) \ge C_1k\log\left(1+\frac{e(p-1)}{k}\right)\,.
\end{equation}
Here, $b_{(1)}, b'_{(1)}$ are the first columns of $B, B'$ with their first components 
replaced by 0.  Thus, for any two distinct matrices $\Sigma$ and $\Sigma' $ belonging to the family
$
\cN'=\{
 \Sigma=I_p+\frac{a}{2}B: \ B \in \cB_1\} 
$ 
we have $
\|\Sigma - \Sigma'\|_1^2 \ge C_4a^2k^2\,
$
for some constant $C_4>0$. Here, $\cN'\subset \cG_0^{(1)}(k)$ thanks to the second inequality in~\eqref{EQ:low:ass}. Also, by Lemma~\ref{LEM:KL}, $
\KL(P_\Sigma, P_{I_p})\le a^2k$ for all $\Sigma\in \cN'$. These remarks and \eqref{EQ:log_card1} imply the suitably modified $(i)$ and $(ii)$ for the choice
$$
a=a_0\left(\frac{1}{n}\log\left(1+\frac{e(p-1)}{k}\right)\right)^{1/2}\,
$$
with $a_0$ small enough. The rest of the proof follows the same lines as the proof of~\eqref{EQ:low1}.
\epr
The lower bound~\eqref{EQ:low2} and Theorem~4 in [CZ] imply that the rate  $R\left((\log p)/n\right)^{\frac{1-q}{2}}$ is optimal on the class $\cG_q^{(1)}(R)$ under the $\|\cdot\|_1$ norm if 
$Rn^{q/2}\le p^{\alpha}$ with some $\alpha<1$. In particular, for $q=0$ this optimality holds under the quite natural condition $k=O( p^{\alpha})$, and no lower bound on $p$ in terms of $n$ is required. Clearly, this is also true when we drop the condition $\Sigma>0$ in the definition of $\cG_q^{(1)}(R)$ and consider a weak $\ell_q$ constraint as in [CZ].

Note that the rate $\psi^{(0)}$ is very similar to the optimal rate in the Gaussian sequence model,  cf. Section 11.5 in~\citet{Joh11}. This is due to the similarity between the vector $\ell_2$ norm and the Frobenius norm. The rate $\psi^{(1)}$ is different but nevertheless has analogous ingredients. Observe also that, in contrast to the remark after Theorem~1 in [CZ], we prove the Frobenius and the $\|\cdot\|_1$-norm lower bounds~\eqref{EQ:low1} and~\eqref{EQ:low2} by exactly the same technique. The key point is the use of the ``$k$-selection lemma" (Lemma A.3 in~\citet{RigTsy11}).  The lower bound~\eqref{EQ:low2} improves upon Theorem~2 in [CZ] in two aspects. First, it does not need the assumption $p\ge n^{\nu}$, $\nu>1$, and provides insight on the presumed optimal rate for any configuration of $n,p,R$. Second, it is established for general loss functions $w$, in particular for the ``in probability" loss that we consider below. The technique used in Theorem~2 of [CZ] is not adapted for this purpose as it applies to special losses derived from $w(t)=t$. 

\medskip

{\bf Approximate sparsity and optimal rates.}
Along with the hard thresholding estimator considered by [CZ], one can use the soft thresholding estimator $\tilde \Sigma$ defined as the matrix with off-diagonal elements 
$$
\tilde \sigma_{ij}=\sign(\sigma_{ij}^*)(|\sigma_{ij}^*|-\tau)_+\,,
$$
where  $\sigma_{ij}^*$ are the elements of the sample covariance matrix $\Sigma^*$, $ \tau>0$ is a threshold, and $(\cdot)_+$ denotes the positive part. The diagonal elements of  $\tilde \Sigma$ are all set to 1 since we consider the classes $\cG_q^{(j)}(R)$, $j=0,1$. Then $\tilde \Sigma=I_p+\tilde \Sigma_{\rm off}$ where $\tilde \Sigma_{\rm off}$ admits the representation (the minimum is taken over all $p\times p$ matrices $S$ with zero diagonal):
$$
\tilde \Sigma_{\rm off} =\argmin_{S: \,{\rm diag}(S)=0} \left\{|S -\Sigma^*|_2^2 + 2\tau |S|_1\right\}\,.
$$
Take the threshold
\begin{equation}\label{tau}
\tau={A\gamma}\sqrt{\frac{\log p}{n}}\,,
\end{equation}
where $A>1$ and $\gamma$ is the constant in the inequality (24) in [CZ]. 
\begin{TH1}\label{TH:up} %
Let $X_1, \ldots, X_n$ be i.i.d. random vectors in $\R^p$ with covariance matrix $\Sigma$ such that (24) in [CZ] holds. 
Assume that $p,n$, and $A$ are such that $\tau \le \delta$, where $\delta$ is the constant introduced after (24) in [CZ]. Then there exists $C_*>0$ such that, with probability at least $1-C_*p^{2-2A^2}$, 
\begin{equation}\label{TH:up:1}
 \big\| \tilde \Sigma -\Sigma\big\|^2\le \min_{S} \left\{ \|S-\Sigma\|^2 + \left(\frac{1+\sqrt{2}}{2}\right)^2A^2\gamma^2 \frac{|S|_0\log p}{n}\right\}\,,
\end{equation}
where $\min_{S}$ denotes the minimum over all $p\times p$ matrices.
\end{TH1}
\noindent {\sc Proof.}  Write $\sigma_{ij}^*=\sigma_{ij}+\xi_{ij}$ where the $\xi_{ij}=\sigma_{ij}^*-\sigma_{ij}$ are zero-mean random variables, $i\ne j$. Thus, considering $\sigma_{ij}^*$ as observations, we have a sequence model in dimension $p(p-1)$. It is easy to see that it is a special case of the trace regression model studied in~\citet{KolLouTsy11} where $A_0$ is a diagonal matrix with the $p(p-1)$ off-diagonal entries of  $\Sigma$ on the diagonal. In the notation of~\citet{KolLouTsy11},  the corresponding matrices $X_i$ are diagonalizations of canonical basis vectors, the norm $\|\cdot\|_{L_2(\Pi)}$ coincides with the  norm $|\cdot|_2$, and ${\rm rank} (B)$ is equal to the number of non-zero entries of diagonal matrix $B$. Thus, Assumption 1 in~\citet{KolLouTsy11} is satisfied with $\mu=1$, and we can apply Theorem~1 in~\citet{KolLouTsy11}. It yields a deterministic statement:
$$
 \big| \tilde \Sigma -\Sigma\big|_2^2\le \min_{S} \left\{ |S-\Sigma|_2^2 + \left(\frac{1+\sqrt{2}}{2}\right)^2\tau^2 |S|_0\right\}
$$
provided $\tau> 2 \max_{i\ne j}|\sigma_{ij}^*-\sigma_{ij}|$. From (24) in [CZ] and a union bound, we obtain that, for $\tau$ defined in \eqref{tau}, this inequality holds with probability greater than $1-C_*p^{2-2A^2}$. 
\epr
 
\begin{cor}
Under the assumptions of Theorem~\ref{TH:up}, for any $0<q<2$, there exist constants $C',C_*>0$ such that with probability at least $1-C_*p^{2-2A^2}$, 
\label{COR:lq_up}
\begin{equation}\label{COR:lq_up:1}
 \big\| \tilde \Sigma - \Sigma \big\|^2\le \min_{S 
 } \left\{ 2\|S-\Sigma\|^2 + C'|S|_q^q \left(\frac{\log p }{n}\right)^{1-\frac{q}{2}}\right\}\,.
\end{equation}
\end{cor}
\noindent {\sc Proof.} Let $|s_{[l]}|, \, l=1, \ldots, p(p-1)$, denote the absolute values of the off-diagonal elements of $S$ ordered in a decreasing order. Note that for any $p\times p$ matrix $S$ and any $0<q<2$ we have $|s_{[l]}|^q \le |S|_q^q/l$. Fix an integer $k\le p(p-1)$. Taking $s'_{ij}=s_{ij}$ if $|s_{ij}|\ge |s_{[k]}|$ and $s'_{ij}=0$ otherwise, we get that for any $S$ there exists a $p\times p$ matrix $S'$ with $|S'|_0=k$ such that
$$
|S-S'|_2^2=\sum_{l >k}s_{[l]}^2 \le |S|_q^2\sum_{l>k}l^{-2/q} \le \frac{|S|_q^2k^{1-2/q}}{2/q-1}\,.
$$
Together with Theorem~\ref{TH:up}, this implies that for any integer $k\le p(p-1)$ we have
$$
\big| \tilde \Sigma -\Sigma\big|_2^2\le  \min_{S} \left\{2|S-\Sigma|_2^2 + \frac{|S|_q^2k^{1-2/q}}{2/q-1} + \left(\frac{1+\sqrt{2}}{2}\right)^2A^2\gamma^2 \frac{k\log p}{n}\right\}\,.
$$
Optimizing the right hand side over $k$ completes the proof. \epr

Note that the oracle inequalities \eqref{TH:up:1} and \eqref{COR:lq_up:1} are satisfied 
for any covariance matrix $\Sigma$, not necessarily for sparse $\Sigma$. They quantify a trade-off between the approximation and sparisty terms.
Their right-hand sides are small if $\Sigma$ is well approximated by a matrix $S$ with a small number of entries or with small $\ell_q$ norm of the off-diagonal elements. If the matrix $\Sigma$ is sparse, $\Sigma\in\cG_q^{(0)}(R)$, the oracle inequalities \eqref{TH:up:1} and \eqref{COR:lq_up:1} imply that 
$$
\sup_{\Sigma\in\cG_q^{(0)}(R)} P_{\Sigma}\left(\| \tilde \Sigma -\Sigma\|>C'' R^{1/2}\left(\frac{\log p}{n}\right)^{\frac1{2}-\frac{q}{4}}\right)\le C_*p^{2-2A^2}\,
$$
for some constant $C''>0$. This also holds when we drop the condition $\Sigma>0$ in the definition of $\cG_q^{(0)}(R)$.
Combining this with Theorem~\ref{TH:low}, we find that the rate  $R^{1/2}\left((\log p)/n\right)^{\frac1{2}-\frac{q}{4}}$ is optimal on the class $\cG_q^{(0)}(R)$ under the Frobenius norm if 
$Rn^{q/2}\le p^{2\alpha}$ with some $\alpha<1$. In particular, for $q=0$ this optimality holds under the condition $k\le p^{2\alpha}$ with some $\alpha<1$.

\vskip 0.2in
\bibliographystyle{ims}
\bibliography{rigolletmain}

\end{document}